\documentclass[12pt]{article}
\usepackage{amsfonts,amscd,a4}
\usepackage{color}
\usepackage{graphicx}
\usepackage{theorem}
\newtheorem{theo}{Theorem}
\newtheorem{prop}[theo]{Proposition}
\newtheorem{lemma}[theo]{Lemma}
\newtheorem{coro}[theo]{Corollary}

{\theorembodyfont{\rm}

}
\newcommand{\bA}{{\bf A}}

\newcommand{\bK}{{\bf K}}

\newcommand{\cA}{{\mathcal A}}

\newcommand{\cF}{{\mathcal F}}
\newcommand{\cG}{{\mathcal G}}

\newcommand{\cJ}{{\mathcal J}}

\newcommand{\cP}{{\mathcal P}}

\newcommand{\cS}{{\mathcal S}}

\newcommand{\eT}{{\sf T}}

\newcommand{\sC}{{\mathbb C}}
\newcommand{\sD}{{\mathbb D}}

\newcommand{\sN}{{\mathbb N}}

\newcommand{\sT}{{\mathbb T}}

\newcommand{\sZ}{{\mathbb Z}}
\newcommand{\qed}{\rule{1ex}{1ex}}

\newcommand{\comm}{\mbox{\rm Comm} \,}

\newcommand{\ess}{{{\rm ess}} \,}

\newcommand{\im}{{\rm im} \,}

\begin{document}
\title{Finite sections of truncated Toeplitz operators}
\author{Steffen Roch}
\date{}
\maketitle
\begin{abstract}
We describe the $C^*$-algebra associated with the finite sections discretization of truncated Toeplitz operators on the model space $K^2_u$ where $u$ is an infinite Blaschke product. As consequences, we get a stability criterion for the finite sections discretization and results on spectral and pseudospectral approximation.

\end{abstract}
{\bf 2010 AMS-MSC:} 65J10, 46L99, 47N40
\section{Truncated Toeplitz operators} \label{ssd4.1}
Let $H^2$ denote the standard Hardy space on the unit disk $\sD$, i.e. the Hilbert space of all holomorphic functions on $\sD$ which have square-summable Taylor coefficients. As usual, we identify $H^2$ with its space of non-tangential boundary functions, which is a closed subspace of the Lebesgue space $L^2(\sT)$ with normalized Lebesgue measure $m$ on the unit circle $\sT$. The orthogonal projection from $L^2(\sT)$ onto $H^2$ is denoted by $P$.

Every function $a \in L^\infty(\sT)$ defines an operator of multiplication on $L^2(\sT)$, which we denote by $aI$. The {\em Toeplitz operator} induced by $a$ is the operator $T(a) := P aI|_{H^2}$, acting from $H^2$ to $H^2$. The Toeplitz operator with generating function $a(z) = z$ is the operator $S$ of forward shift, $(Sf)(z) = zf(z)$. Its adjoint $S^*$, the backward shift operator, is given by $(S^*f)(z) = z^{-1} (f(z) - f(0))$.

Let $u$ be a non-constant inner function, i.e., $u$ is holomorphic on $\sD$ and $|u(t)| = 1$ for $t \in \sT$ (the following becomes trivial when $u$ is constant). The subspace $K_u^2 := H^2 \ominus u H^2$ is a proper nontrivial invariant subspace of $S^*$. Conversely, every proper nontrivial invariant subspace of $S^*$ is of this form by a celebrated theorem of Beurling. The spaces $K^2_u$ are also known as {\em model spaces}. We denote the orthogonal projection from $L^2(\sT)$ onto $K^2_u$ by $P_u$. If $M_u$ and $M_{\bar{u}}$ denote the operators of multiplication by $u$ and $\bar{u}$ on $L^2(\sT)$, then $P_u = P - M_u P M_{\bar{u}}$.

For $a \in L^\infty(\sT)$, the {\em truncated Toeplitz operator} (TTO for short) generated by $a$ is the operator $T_u(a) := P_u aI|_{K^2_u}$ acting from $K^2_u$ to $K^2_u$. Truncated Toeplitz operators share many of their properties with their relatives, the Toeplitz operators on $H^2$, to which we here  sometimes refer as {\em classical} Toeplitz operators, but there are also some striking differences. For example, the function $a$ is in general not uniquely determined by the operator $T_u(a)$ it generates, and the truncated shift $S_u$, i.e., the TTO with generating function $a(z) = z$, is the sum of a unitary and a compact operator (hence a Fredholm operator with index $0$), whereas its classical counterpart $S$ is a proper isometry (and a Fredholm operator of index $-1$). Moreover, whereas the spectrum of the classical shift $S$ is the closed unit disk $\overline{\sD}$, the spectrum of the truncated shift $S_u$ coincides with the so-called {\em spectrum}
\[
\sigma(u) := \{ \lambda \in \overline{\sD} : \liminf_{z \to \lambda} |u(z)| = 0 \}
\]
of the inner function $u$ (\cite[Lemma 2.5]{Sar3}).

The following collection of positive results is taken from and proved in \cite{GRW1}. We denote by $\eT_u(C)$ the smallest closed $C^*$-subalgebra of $L(K^2_u)$ which contains the truncated shift $S_u$ and the identity operator (the notation will be justified by assertion $(d)$ in the theorem below). Further, we write $\comm \cA$ for the commutator ideal of a $C^*$-algebra $\cA$, i.e., for the smallest closed ideal of $\cA$ which contains all commutators $ab-ba$ with $a, \, b \in \cA$. The essential spectrum of an operator $A$ is denoted by $\sigma_\ess (A)$, and its essential norm by $\|A\|_\ess$.
\begin{theo} \label{td4.1}
Let $u$ be a non-constant inner function. Then \\[1mm]
$(a)$ for $a, \, b \in C(\sT)$, $T_u(a) T_u(b) - T_u(ab)$ is compact. \\[1mm]
$(b)$ $\comm (\eT_u(C)) = K(K^2_u)$. \\[1mm]
$(c)$ $\eT_u(C)/K(K^2_u)$ is $^*$-isomorphic to $C(\sigma(u) \cap \sT)$. \\[1mm]
$(d)$ for $a \in C(\sT)$, the TTO $T_u(a)$ is compact if and only if $a(\sigma(u) \cap \sT) = \{0\}$. \\[1mm]
$(e)$ $\eT_u(C) = \{T_u(a) + K : a \in C(\sT), \, K \in K(K^2_u)\}$. \\[1mm]
$(f)$ for $a \in C(\sT)$, $\sigma_\ess (T_u(a)) = a(\sigma_\ess (S_u))$. \\[1mm]
$(g)$ for $a \in C(\sT)$, $\|T_u(a)\|_\ess = \sup \{ |a(t)| : t \in \sigma(u) \cap \sT\}$. \\[1mm]
$(h)$ Every operator in $\eT_u(C)$ is the sum of a normal and a compact operator. \\[1mm]
Moreover,
\[
\{0\} \longrightarrow K(K^2_u) \stackrel{\rm id}{\longrightarrow} \eT_u(C) \stackrel{\pi}{\longrightarrow} C(\sigma(u) \cap \sT) \longrightarrow \{0\}
\]
is a short exact sequence, with the mapping $\pi$ given by $T_u(a) + K \mapsto a|_{\sigma(u) \cap \sT}$.
\end{theo}
Here is an outline of the contents of the paper. In Section \ref{ssd4.2} we will single out a sequence $(P_{u_n})$ of finite rank projections which converge strongly to the identity operator on $K^2_u$. The operator $P_{u_n} T_u(a) P_{u_n}$ is considered as a finite section of the truncated Toeplitz operator $T_u(a)$. In Section \ref{ssd4.4}, Theorem \ref{td4.14}, we describe the $C^*$-algebra $\cS(\eT_u(C))$ generated by all sequences of the form $(P_{u_n} T_u(a) P_{u_n})_{n \ge 1}$ with $a$ a continuous function. This description is based of a formula of Widom-type that we will derive in Section \ref{ssd4.3}. As consequences of Theorem \ref{td4.14}, we get a stability criterion and results on spectral approximation. The stability criterion (Theorem \ref{td4.15}) says that a sequence $(A_n)$ in $\cS(\eT_u(C))$ is stable if and only if its strong limit $A$ is invertible. In the case when the $A_n$ are the finite sections of a truncated Toeplitz operator, this result is due to Treil \cite{Trl1}. One advantage of Theorem \ref{td4.15} is that it implies (without any additional effort) results on spectral and pseudospectral approximation as well as on the asymptotic behavior of the small singular values of $A_n$; see the end of Section \ref{ssd4.4}.
\section{A filtration and Widom's identity} \label{ssd4.2}
Recall that a {\em filtration} on a Hilbert space $H$ is a sequence $\cP = (P_n)$ of orthogonal projections of finite rank on $H$ which converges strongly to the identity operator on $H$. To define a filtration on the model space $K^2_u$ we specify $u$ to be a Blaschke product, as follows. A single {\em Blaschke factor} is a function on the unit disk of the form
\begin{equation} \label{ed4.1a}
b_\lambda (z) := \left\{
\begin{array}{cll}
z & \mbox{if} & \lambda = 0, \\
\frac{\lambda - z}{1 - \bar{\lambda}z} \frac{|\lambda|}{\lambda} & \mbox{if} & \lambda \in \sD \setminus \{0\}.
\end{array}
\right.
\end{equation}
A {\em Blaschke product} is then a function
\begin{equation} \label{ed4.2}
u = \prod_{\lambda \in \sD \, : \, k(\lambda) > 0} b_\lambda^{k(\lambda)}
\end{equation}
which satisfies the {\em Blaschke condition}
\begin{equation} \label{ed4.3}
\sum_{\lambda \in \sD} k(\lambda) (1 - |\lambda|) < \infty.
\end{equation}
If $u$ is a finite Blaschke product, i.e., if $u$ is of the form (\ref{ed4.2}) with $k(\lambda) = 0$ for all but finitely many $\lambda \in \sD$, then (\ref{ed4.3}) is satisfied. Conversely, if (\ref{ed4.3}) holds, then every disk $\{ z \in \sD : |z| \le r\}$ with $0 < r < 1$ contains only finitely many $\lambda$ with $k(\lambda) \neq 0$. Thus, if $u$ in (\ref{ed4.2}) is an (infinite) Blaschke product, the number of its non-one factors is countable. We order the $\lambda$ with $k(\lambda) \neq 0$ in a sequence $(\lambda_k)_{k \ge 1}$ in such a way that $|\lambda_k| \le |\lambda_{k+1}|$ for all $k$. Then  (\ref{ed4.2}) and (\ref{ed4.3}) can be written as
\begin{equation} \label{ed4.4}
u = \prod_{k=1}^\infty b_{\lambda_k} \quad \mbox{with} \quad \sum_{k=1}^\infty (1 - |\lambda_k|) < \infty.
\end{equation}
Blaschke products are inner functions. If $u$ is an infinite Blaschke product, we use its product representation (\ref{ed4.4}) to define a filtration on the associated model space $K^2_u$. For $n \ge 1$, set $u_n := \prod_{k=1}^n b_{\lambda_k}$ and let $P_{u_n}$ denote the orthogonal projection from $L^2(\sT)$ onto $K^2_{u_n}$. The projections $P_{u_n}$ own the following properties.
\begin{prop} \label{pd4.5}
$(a)$ The projections $P_{u_n}$ have a finite spatial rank. \\[1mm]
$(b)$ $P_{u_n} \to P_u$ on $L^2(\sT)$ strongly as $n \to \infty$. \\[1mm]
$(c)$ $P_{u_m}P_{u_n} = P_{u_n}P_{u_m} = P_{u_{\min \{m, \, n\}}}$ for $m, \, n \ge 1$. \\[1mm]
$(d)$ $P_{u_n} P_u = P_u P_{u_n} = P_{u_n}$ and $P_{u_n} P = P P_{u_n} = P_{u_n}$ for $n \ge 1$.
\end{prop}
{\bf Hints to the proof.} Assertion $(a)$ is the "Lemma on Finite Dimensional Subspaces" in \cite[p. 33]{Nik1}. For assertion $(b)$ observe that the $u_n$ converge to $u$ uniformly on compact subsets of $\sD$ by the "Lemma on Blaschke Products" in \cite[p. 280]{Nik1}. Thus, the model spaces $K^2_{u_n}$ converge to $K^2_u$ by the "Theorem on Lower Limits" in \cite[p. 34]{Nik1} (note that the limit of the $u_n$ exists; so the set of its limit points contains only one element, which clearly coincides with the greatest common divisor in the formulation of that theorem). Thus, $P_{u_n} \to P_u$ on $L^2(\sT)$ strongly by the definition in \cite[p. 34]{Nik1}.

Finally, $u_n$ is a divisor of $u$ and of $u_m$ for $m \ge n$. Corollary 8 in \cite[p. 19]{Nik1} then implies that $K^2_{u_n} \subseteq K^2_{u_m} \subseteq K^2_u \subseteq H^2$, whence assertions $(c)$ and $(d)$. \hfill \qed \\[3mm]
Thus, the restrictions of the projections $P_{u_n}$ to the model space $K^2_u$ form a filtration on $K^2_u$ by the preceding proposition. We denote this filtration by $\cP_u$ and write $\cP_u = (P_{u_n})_{n \ge 1}$, not distinguishing between a projection $P_{u_n}$ on $L^2(\sT)$ and its restriction to $K^2_u$.

The study of the finite sections discretization (FSD for short) for (classical) Toeplitz operators is dominated by Widom's identity
\begin{equation} \label{e14.20aa}
P_n T(ab) P_n = P_n T(a) P_n T(b) P_n + P_n H(a)H(\tilde{b}) P_n + R_n H(\tilde{a})H(b) R_n.
\end{equation}
To explain this identity, we need some notation. Let $J : L^2(\sT) \to L^2(\sT)$ denote the operator $(Jf)(t) = t^{-1} f(t^{-1})$. One easily checks that, for every function $a \in L^\infty(\sT)$, $JaJ$ is the operator of multiplication by the function $\tilde{a}(t) := a(t^{-1})$. Then $H(a) := P aI J|_{H^2}$ is the (classical) Hankel operator, $R_n := H(t^n)$, and $P_n := R_n^2$. Note that $P_n$ is just the orthogonal projection onto the linear span of the functions $t^n$ with $n \in \{0, \, 1, \, \ldots, \, n-1\}$ and that $R_n = R_n^*$.

Our goal is to achieve a comparable identity for the finite sections of TTO, where the role of the $R_n$ in Widom's identity is played by the operators
\begin{equation} \label{ed4.6}
R_{u_n} := P_{u_n} M_{u_n} J \quad \mbox{and} \quad  R_{u_n}^* = J M_{\overline{u_n}} P_{u_n}.
\end{equation}
\begin{theo}[Widom's identity for TTO] \label{td4.7}
Let $a, \, b \in L^\infty(\sT)$. Then
\begin{eqnarray*}
\lefteqn{P_{u_n} T_u(ab) P_{u_n}} \\
&& = P_{u_n} T_u(a) P_{u_n} T_u(b) P_{u_n} + P_{u_n} H(a) H(\tilde{b}) P_{u_n} + R_{u_n} H(\tilde{a}) H(b) R_{u_n}^*.
\end{eqnarray*}
\end{theo}
{\bf Proof.} By Proposition \ref{pd4.5} $(d)$ and since $P_{u_n} = P - M_{u_n} P M_{\overline{u_n}}$, we find
\begin{eqnarray*}
\lefteqn{P_{u_n} T_u(a) P_{u_n} T_u(b) P_{u_n}} \\
&& = P_{u_n} P a P_{u_n} b P P_{u_n} \\
&& = P_{u_n} P a P b P P_{u_n} - P_{u_n} a M_{u_n} P M_{\overline{u_n}} b P_{u_n} \\
&& = P_{u_n} P ab P P_{u_n} - P_{u_n} P a Q b P P_{u_n} - P_{u_n} M_{u_n} a P b M_{\overline{u_n}} P_{u_n} \quad \mbox{with} \; Q := I - P \\
&& = P_{u_n} ab P_{u_n} - P_{u_n} P a Q J^2 Q b P P_{u_n} \\
&& \qquad - \; P_{u_n} M_{u_n} P a P b M_{\overline{u_n}} P_{u_n} - P_{u_n} M_{u_n} Q a P b M_{\overline{u_n}} P_{u_n} \\
&& = P_{u_n} T_u(ab) P_{u_n} - P_{u_n} H(a) H(\tilde{b}) P_{u_n} - P_{u_n} M_{u_n} Q a P b M_{\overline{u_n}} P_{u_n}.
\end{eqnarray*}
In the last line we used that $P_{u_n} M_{u_n} Pf \in P_{u_n} (u_n H^2) = 0$ for $f \in L^2(\sT)$, whence $P_{u_n} M_{u_n} P = 0$. Then $P M_{\overline{u_n}} P_{u_n} = (P_{u_n} M_{u_n} P )^* = 0$, too, and we can proceed with
\begin{eqnarray*}
\lefteqn{P_{u_n} T_u(a) P_{u_n} T_u(b) P_{u_n}} \\
&& = P_{u_n} T_u(ab) P_{u_n} - P_{u_n} H(a) H(\tilde{b}) P_{u_n} \\
&& \qquad - \; P_{u_n} M_{u_n} Q a P b P M_{\overline{u_n}} P_{u_n} - P_{u_n} M_{u_n} Q a P b Q M_{\overline{u_n}} P_{u_n} \\
&& = P_{u_n} T_u(ab) P_{u_n} - P_{u_n} H(a) H(\tilde{b}) P_{u_n} - P_{u_n} M_{u_n} Q a P b Q M_{\overline{u_n}} P_{u_n} \\
&& = P_{u_n} T_u(ab) P_{u_n} - P_{u_n} H(a) H(\tilde{b}) P_{u_n} - P_{u_n} M_{u_n} J^2 Q a P b Q J^2 M_{\overline{u_n}} P_{u_n} \\
&& = P_{u_n} T_u(ab) P_{u_n} - P_{u_n} H(a) H(\tilde{b}) P_{u_n} - R_{u_n} J Q a P b Q J R_{u_n}^*.
\end{eqnarray*}
Since $J Q a P b Q J = H(\tilde{a}) H(b)$, this is the assertion  \hfill \qed
\section{Hankel operators by Blaschke products} \label{ssd4.3}
The operators $R_n$ in Widom's identity (\ref{e14.20aa}) can be identified with the (classical) Hankel operators $H(t^n)$. Similarly,
\begin{eqnarray*}
R_{u_n} & = & P_{u_n} M_{u_n} J \\
& = & (P - M_{u_n} P M_{\overline{u_n}}) M_{u_n} J \\
& = & (P M_{u_n} - M_{u_n} P) J \\
& = & (P M_{u_n} - P M_{u_n} P) J \quad \mbox{(since $M_{u_n} Pf \in H^2$ for $f \in H^2$)} \\
& = & P M_{u_n} Q J,
\end{eqnarray*}
which can be identified with the (classical) Hankel operator $H(M_{u_n})$ on $H^2$. Analogously, $R_{u_n}^*$ can be identified with $H(M_{u_n})^* = H(\overline{\widetilde{M_{u_n}}})$. We will see that the operators $R_{u_n}$ in Widom's identity for TTO play a quite different role compered with the $R_n$ in (\ref{e14.20aa}). We start with some general properties of Hankel operators generated by inner functions.
\begin{prop} \label{pd4.8}
Let $u \in H^\infty$ and $|u(t)| = 1$ for $t \in \sT$. Then the Hankel operator $H(u) = PuQJ$ is a partial isometry, the range and initial projection of which are given by $H(u) H(u)^* = P - u P \bar{u} I$ and $H(u)^* H(u) = P - \bar{\tilde{u}} P \tilde{u}I$.
\end{prop}
{\bf Proof.} Since $u\bar{u} = 1$ and $PuP = uP$,
\begin{eqnarray} \label{ed4.9}
\lefteqn{H(u) H(u)^* = PuQJ JQ\bar{u} P} \nonumber \\
&& = PuQ\bar{u}P = P u\bar{u} P - PuP\bar{u}P = P - uP\bar{u}I.
\end{eqnarray}
Using this identity and $P\bar{u}Q = 0$ we obtain
\begin{eqnarray*}
\lefteqn{H(u) H(u)^* H(u) = PuQJ - u P\bar{u}PuQJ} \\
&&  = PuQJ - u P\bar{u}uQJ + u P\bar{u}QPuQJ = PuQJ = H(u),
\end{eqnarray*}
i.e., $H(u)$ is an isometry. Replacing $u$ in (\ref{ed4.9}) by $\bar{\tilde{u}}$ (which is also in $H^\infty$) and taking into account that $H(\bar{\tilde{u}}) = H(u)^*$, the identity for $H(u)^* H(u)$ follows. \hfill \qed
\begin{coro} \label{cd4.10}
$R_{u_n} R_{u_n}^* = P_{u_n}$, $R_{u_n}^* R_{u_n} = P_{\overline{\widetilde{u_n}}}$, $P_{u_n} R_{u_n} = R_{u_n}$, $R_{u_n}^* P_{u_n} = R_{u_n}^*$.
\end{coro}
The following convergence result for the $R_{u_n}$ is in sharp contrast with the $R_n$ in Widom's identity (\ref{e14.20aa}), which converge weakly to zero.
\begin{theo} \label{td4.11}
$R_{u_n} = H(M_{u_n}) \to H(M_u)$ $^*$-strongly.
\end{theo}
In the proof of this result, we will make use of the following well known assertion.
\begin{lemma} \label{l1}
Let $A_n, \, A$ be bounded linear operators on a Hilbert space $H$. If $A_n \to A$ weakly and $\|A_nx\| \to \|Ax\|$ for all $x \in H$, then $A_n \to A$ strongly.
\end{lemma}
{\bf Proof.}  It is clearly sufficient to prove the following fact for elements $x_n, \, x$ of $H$: if $x_n \to x$ weakly and $\|x_n\| \to \|x\|$, then $\|x_n - x\| \to 0$. This follows from
\begin{eqnarray*}
\|x_n - x\|^2 & = & \langle x_n - x, \, x_n - x \rangle \\
& = & \langle x_n, \, x_n \rangle - \langle x_n, \, x \rangle - \langle x, \, x_n \rangle + \langle x, \, x \rangle \\
& = & \|x_n\|^2 + \|x\|^2 - \langle x_n, \, x \rangle - \langle x, \, x_n \rangle
\end{eqnarray*}
which goes to 0 by hypothesis. \hfill \qed \\[3mm]
{\bf Proof of Theorem \ref{td4.11}.} We first show that $R_{u_n} = H(M_{u_n}) \to H(M_u)$ weakly. Indeed, the uniform convergence of $u_n$ to $u$ on compact subsets of $\sD$ implies the convergence of the $k$th Taylor coefficient of $u_n$ to the $k$th Taylor coefficient of $u$ for every $k \in \sZ^+$. Together with the uniform boundedness of the operators $H_{u_n}$, this fact implies the weak convergence of $H(M_{u_n})$ to $H(M_u)$.

Next we show that $\|H(M_{u_n}) x\| \to \|H(M_u) x\|$ for every $x \in H^2$. Once this is done, the strong convergence of $H(M_{u_n})$ to $H(M_u)$ follows from Lemma \ref{l1}.

We start with showing that
\begin{equation} \label{ed4.12}
P_{\overline{\widetilde{u_n}}} \to P_{\bar{\tilde{u}}} \quad \mbox{strongly as} \; n \to \infty.
\end{equation}
Indeed, let $b_\lambda$ be a single Blaschke factor as in (\ref{ed4.1a}). For $t \in \sT$ and $\lambda \neq 0$ we then have
\[
\overline{\widetilde{b_\lambda}} (t)
= \overline{\frac{\lambda - t^{-1}}{1 - \bar{\lambda}t^{-1}} \frac{|\lambda|}{\lambda}}
= \frac{\bar{\lambda} - t}{1 - \lambda t} \frac{|\lambda|}{\bar{\lambda}} = b_{\bar{\lambda}} (t).
\]
For $\lambda = 0$, the equality $\overline{\widetilde{b_\lambda}} = b_{\bar{\lambda}}$ on $\sT$ is evident. Moreover, if $(\lambda_n)$ is a sequence in $\sD$ satisfying the Blaschke condition, then the sequence $(\overline{\lambda_n})$ also satisfies this condition. So we can apply the assertion of Proposition \ref{pd4.5} $(b)$ to the functions
$\prod_{k=1}^\infty b_{\overline{\lambda_k}}$ and $\prod_{k=1}^n b_{\overline{\lambda_k}}$ in place of $u$ and $u_n$ to get the assertion (\ref{ed4.12}). From (\ref{ed4.12}) and Corollary \ref{cd4.10} we then conclude that
\[
H(M_{u_n})^* H(M_{u_n}) \to H(M_u)^* H(M_u) \quad \mbox{strongly},
\]
from which we obtain
\[
\langle H(M_{u_n})^* H(M_{u_n}) x, \, x \rangle - \langle H(M_u)^* H(M_u) x, \, x \rangle = \|H(M_{u_n}) x\|^2 - \|H(M_u) x\|^2 \to 0
\]
for every $x \in H^2$. This proves the strong convergence of $R_{u_n} = H(M_{u_n})$ to $H(M_u)$. The strong convergence of the adjoint operators follows as above, by working with the Blaschke product $\bar{\tilde{u}}$ in place of $u$. \hfill \qed
\begin{coro} \label{cd4.13}
Let $L$ be a compact operator on $H^2$. Then
\[
\|R_{u_n} L R_{u_n}^* - P_{u_n} H(M_u) L H(M_u)^* P_{u_n}\| \to 0 \quad \mbox{as} \; n \to \infty.
\]
\end{coro}
Indeed, since $L$ and $H(M_u) L H(M_u)^*$ are compact, we derive from Proposition \ref{pd4.5} $(b)$ and Theorem \ref{td4.11} that both sequences $(R_{u_n} L R_{u_n}^*)$ and $(P_{u_n} H(M_u) L H(M_u)^* P_{u_n})$ converge to $H(M_u) L H(M_u)^*$ in the norm. \hfill \qed
\section{The algebra of the FSD for TTO} \label{ssd4.4}
With every filtration $\cP = (P_n)$ on a Hilbert space $H$, there are naturally associated some algebraic objects. By $\cF^\cP$ we denote the set of all sequences $\bA = (A_n)$ of operators $A_n : \im P_n \to \im P_n$ for which the sequence $(A_n P_n)$ converges $^*$-strongly to some operator $W(\bA)$ on $H$. Provided with element-wise defined operations and the supremum norm, $\cF^\cP$ becomes a $C^*$-algebra, the set $\cG^\cP$ of all sequences in $\cF^\cP$ which converge to 0 in the norm is a closed ideal of $\cF^\cP$, and the mapping $W$, also called the consistency map of the filtration $\cP$, is a $^*$-homomorphism from $\cF^\cP$ to the algebra $L(H)$ of the bounded linear operators on $H$.

We prepare the proof of Theorem \ref{td4.14} below by an assertion of independent interest.
\begin{prop} \label{p12.20}
Let $\cP = (P_n)$ be a filtration on a Hilbert space $H$. Then the ideal $\cG^\cP$ is contained in the smallest closed subalgebra $\cJ$ of $\cF^\cP$ which contains all sequences $(P_n K P_n)$ with $K$ compact if and only if $\cP$ is injective.
\end{prop}
{\bf Proof.} The "only if"-part of the assertion is evident. For the "if"-part we are going to show that, for each $n_0 \in \sN$, there is a sequence $(G_n)$ in $\cJ$ with $G_{n_0} \neq 0$ and $G_n = 0$ for all $n \neq n_0$. Since the matrix algebras $\sC^{k \times k}$ are simple, this fact already implies that each sequence $(G_n)$ with arbitrarily prescribed $G_{n_0} \in L(\im P_{n_0})$ and $G_n = 0$ for $n \neq n_0$ belongs to $\cJ$. Since  $\cG^\cP$ is generated by sequences of this special form, the assertion follows.

For $n_0 \in \sN$, put
\[
\sN_< := \{ n \in \sN : \im P_n \cap \im P_{n_0} \; \mbox{is a proper subspace of} \; \im P_{n_0} \},
\]
and set $\sN_> := \sN \setminus (\{n_0\} \cup \sN_<)$. The set $\sN_<$ is at most countable, and none of the closed linear spaces $\im P_n \cap \im P_{n_0}$ has interior points relative to $\im P_{n_0}$. By the Baire category theorem,
$\cup_{n \in \sN_<} (\im P_n \cap \im P_{n_0})$ is a proper subset of $\im P_{n_0}$. Choose a unit vector
\[
f \in \im P_{n_0} \setminus \cup_{n \in \sN_<} (\im P_n \cap \im P_{n_0}).
\]
Then $\|P_n f\| < 1$ for all $n \in \sN_<$ by the Pythagoras theorem. (Indeed, otherwise $\|P_n f\| = 1$, and the equality $1 = \|f\|^2 = \|P_n f\|^2 + \|f - P_n f\|^2$ implies $f = P_n f$, whence $f \in \im P_n$.)

Let $Q_n:= I - P_n$. If $n \in \sN_>$, then $\im P_n \cap \im P_{n_0} = \im P_{n_0}$ by the definition of $\sN_>$. Thus, $\im P_{n_0} \subseteq \im P_n$, and since no two of the projections $P_n$ coincide, this implies that $\im P_{n_0}$ is a proper subspace of $\im P_n$ and $\im Q_n$ is a proper subspace of $\im Q_{n_0}$ for $n \in \sN_>$. Again by the Baire category theorem, $\cup_{n \in \sN_>} \im Q_n$ is a proper subset of $\im Q_{n_0}$. Choose a unit vector
\[
g \in \im Q_{n_0} \setminus \cup_{n \in \sN_>} \im Q_n.
\]
As above,  $\|Q_n g\| < 1$ for all $n \in \sN_>$. Consider the operator $K : x \mapsto \langle x, \, g \rangle f$ on $H$. Its adjoint is the operator $x \mapsto \langle x, \, f \rangle g$, and
\[
P_n K Q_n K^* P_n x = \langle P_n x, \, f \rangle  \, \langle Q_n g, \, g \rangle \, P_n f
= \langle x, \, P_n f \rangle \, \|Q_n g\|^2 P_n f.
\]
If $n \in \sN_<$ then $\|P_n f\| < 1$, and if $n \in \sN_>$ then $\|Q_n g\| < 1$ by construction. In both cases, $\|P_n K Q_n K^* P_n\| < 1$. In case $n = n_0$, the operator $P_n K Q_n K^* P_n x = \langle x, \, f \rangle f$ is an orthogonal projection of norm 1, which we call $P$. The sequence $\bK := (P_n K Q_n K^* P_n)$ belongs to the algebra $\cJ$ since
\[
(P_n K Q_n K^* P_n) = (P_n KK^*P_n) - (P_n K P_n) \, (P_n K^* P_n).
\]
As $r \to \infty$, the powers $\bK^r$ converge in the norm of $\cF^\cP$ to the sequence $(G_n)$ with $G_{n_0} = P \neq 0$ and $G_n = 0$ if $n \neq n_0$. Indeed, since $P_n \to I$ strongly, one has $\|Q_n g\| < 1/2$ for $n$ large enough, whence $\|P_n K Q_n K^* P_n\| < 1/2$ for these $n$, whereas $\|P_n K Q_n K^* P_n\| < 1$ for the remaining (finitely many) $n$ as seen above. Since $\bK^r \in \cJ$ and $\cJ$ is closed, the sequence $(G_n)$ has the claimed properties. \hfill \qed \\[3mm]
The goal in this section is to study the FSD of TTO with respect to the filtration $\cP_u$. In accordance with Theorem \ref{td4.1} $(e)$, we define the corresponding (full) algebra of the FSD as the smallest closed subalgebra $\cS(\eT_u(C))$ of $\cF^{\cP_u}$ which contains all sequences $(P_{u_n} (T_u(a) + K)P_{u_n})_{n \ge 1}$ with $a \in C(\sT)$ and $K \in K(K^2_u)$.
\begin{theo} \label{td4.14}
$\cS(\eT_u(C))$ consists of all sequences $(P_{u_n} (T_u(a) + K) P_{u_n} + G_n)$ with $a \in C(\sT), \, K \in K(K^2_u)$ and $(G_n) \in \cG^{\cP_u}$.
\end{theo}
{\bf Proof.} The proof runs parallel to that of Theorem \cite[1.53]{HRS2}; so we address to some main steps only.

For a moment, let $\cS_1$ denote the set of all sequences of the mentioned form. The sequences $(P_{u_n} (T_u(a) + K) P_{u_n})$ are contained in $\cS(\eT_u(C))$ by definition, and since the filtration $\cP_u$ is injective, we conclude from Proposition \ref{p12.20} that the ideal $\cG^{\cP_u}$ of the zero sequences is also contained in $\cS(\eT_u(C))$. Thus, $\cS_1 \subseteq \cS(\eT_u(C))$.

For the reverse inclusion, we prove that $\cS_1$ is a closed subalgebra of $\cS(\eT_u(C))$. If $a, \, b \in C(\sT)$ then, by Theorem \ref{td4.7} (Widom's identity) and Corollary \ref{cd4.13},
\begin{eqnarray*}
\lefteqn{P_{u_n} T_u(a) P_{u_n} \cdot P_{u_n} T_u(b) P_{u_n}} \\
&& = P_{u_n} T_u(ab) P_{u_n} - P_{u_n} H(a) H(\tilde{b}) P_{u_n} - R_{u_n} H(\tilde{a}) H(b) R_{u_n}^* \\
&& = P_{u_n} T_u(ab) P_{u_n} - P_{u_n} (H(a) H(\tilde{b}) + H(M_u) H(\tilde{a}) H(b) H(M_u)^*)P_{u_n} + G_n \\
&& = P_{u_n} T_u(ab) P_{u_n} - P_{u_n} K P_{u_n} + G_n
\end{eqnarray*}
with a compact operator $K$ and a sequence $(G_n) \in \cG^{\cP_u}$. Thus,
\[
(P_{u_n} T_u(a) P_{u_n}) \, (P_{u_n} T_u(b) P_{u_n}) \in \cS_1.
\]
It follows now in a standard way that $\cS_1$ is an algebra. To prove that $\cS_1$ is closed, let $((P_{u_n} (T_u(a_m) + K_m) P_{u_n} + G_n^m)_{n \ge 1})_{m \ge 1}$ be a sequence in $\cS_1$ which converges in the norm of $\cF^{\cP_u}$. Let $W$ denote the consistency map of the filtration $\cP_u$, i.e., $W((A_n)) = \mbox{s-lim} \, A_n P_{u_n}$. Then
\[
(W((P_{u_n} (T_u(a_m) + K_m) P_{u_n} + G_n^m)_{n \ge 1}))_{m \ge 1} = (T_u(a_m) + K_m)_{m \ge 1}
\]
is a Cauchy sequence in $\eT_u(C)$. Since $\eT_u(C)$ is a closed algebra, this sequence converges in $\eT_u(C)$. The limit of this sequence is of the form $T_u(a) + K$ with $a \in C(\sT)$ and $K$ compact by Theorem \ref{td4.1} $(e)$. (But note that the representation of the limit in that form is not unique by Theorem \ref{td4.1} $(d)$.) It is now easy to see that
\[
(P_{u_n} (T_u(a_m) + K_m) P_{u_n})_{n \ge 1} \to (P_{u_n} (T_u(a) + K) P_{u_n})_{n \ge 1}
\]
in the norm of $\cF^{\cP_u}$ as $m \to \infty$. Then, finally, the sequence $((G_n^m)_{n \ge 1})_{m \ge 1}$ converges; its limit is in $\cG^{\cP_u}$. \hfill \qed
\begin{theo} \label{td4.15}
A sequence $\bA = (A_n) \in \cS(\eT_u(C))$ is stable if and only if the operator $W(\bA) = \mbox{\rm s-lim} \, A_n P_{u_n}$ is invertible.
\end{theo}
{\bf Proof.} By Theorem \ref{td4.14}, we have to show that the sequence $\bA := (P_{u_n} (T_u(a) + K) P_{u_n})$ (with $a \in C(\sT)$ and $K \in K(K^2_u)$) is stable if and only if the operator $A := T_u(a) + K$ is invertible. Since the stability of $\bA$ implies the invertibility of $A$ by Polski's theorem, we are left with the reverse implication.

So let $A$ be invertible. By inverse closedness of $C^*$-algebras, $A^{-1} \in \eT_u(C)$; hence, $A^{-1} = T_u(b) + L$ with a certain function $b \in C(\sT)$ and a compact operator $L$ by Theorem \ref{td4.1} $(e)$. Using Widom's identity as in the proof of the previous theorem and employing assertion $(a)$ of Theorem \ref{td4.1}, we conclude that there are compact operators $R_1, \, R_2$ and a sequence $(G_n) \in \cG^{\cP_u}$ such that
\begin{eqnarray*}
\lefteqn{P_{u_n} T_u(a) P_{u_n} \cdot P_{u_n} T_u(b) P_{u_n} - P_{u_n}} \\
&& = P_{u_n} (T_u(ab) - I) P_{u_n} - P_{u_n} R_1 P_{u_n} + G_n \\
&& = P_{u_n} R_2 P_{u_n} - P_{u_n} R_1 P_{u_n} + G_n.
\end{eqnarray*}
Thus, the sequence $(P_{u_n} T_u(b) P_{u_n})$ is a right inverse of the sequence $(P_{u_n} T_u(a) P_{u_n})$ modulo the ideal
\[
\cJ := \{(P_{u_n} K P_{u_n} + G_n) : K \in K(K^2_u), \, (G_n) \in \cG^{\cP_u}\}
\]
of $\cS(\eT_u(C))$. A similar computation shows that it is a left inverse modulo $\cJ$, too. Then  $(P_{u_n} T_u(b) P_{u_n})$ is also an inverse of $(P_{u_n} A P_{u_n})$ modulo $\cJ$. Now the assertion follows from the Lifting Theorem \cite[5.37]{HRS2} in its simplest form, i.e., with $\cJ/\cG^{\cP_u}$ consisting of one elementary ideal only. \hfill \qed \\[3mm]
The following is certainly the most important consequence of Theorem \ref{td4.15}. The definition of a fractal algebra is in \cite{HRS2}.
\begin{coro} \label{cd4.16}
The algebra $\cS(\eT_u(C))$ is fractal.
\end{coro}
Note that the consistency map of $\cP_u$ is fractal; so the assertion of the corollary follows from Theorems \cite[1.69]{HRS2} and \ref{td4.15}.  \hfill \qed \\[3mm]
Sequences in fractal algebras are distinguished by their excellent convergence properties. To mention only a few of them, let $\sigma (a)$ denote the spectrum of an element $a$ of a $C^*$-algebra with identity element $e$, write $\sigma_2(a)$ for the set of the singular values of $a$, i.e., $\sigma_2(a)$ is the set of all non-negative square roots of elements in the spectrum of $a^*a$ and finally, for $\varepsilon > 0$, let $\sigma^{(\varepsilon)} (a)$ refer to the $\varepsilon$-pseudospectrum of $a$, i.e. to the set of all $\lambda \in \sC$ for which $a - \lambda e$ is not invertible or $\|(a - \lambda e)^{-1}\| \ge 1/\varepsilon$. Let further
\[
d_H(M, \, N) := \max \, \{ \max_{m \in M} \min_{n \in N} |m-n|, \, \max_{n \in N} \min_{m \in M} |m-n| \}
\]
denote the Hausdorff distance between the non-empty compact subsets $M$ and $N$ of the complex plane.
\begin{theo} \label{t95.41}
Let $(A_n)$ be a sequence in $\cS(\eT_u(C))$ with strong limit $A$. Then the following set-sequences converge with respect to the Hausdorff distance as $n \to \infty \!:$ \\[1mm]
$(a)$ $\sigma (A_n) \to \sigma (A)$ if $(A_n)$ is self-adjoint; \\[1mm]
$(b)$ $\sigma_2 (A_n) \to \sigma_2 (A)$; \\[1mm]
$(c)$ $\sigma^{(\varepsilon)} (A_n) \to \sigma^{(\varepsilon)} (A)$.
\end{theo}
The proof follows immediately from the stability criterion in Theorem \ref{td4.15} and from Theorems 3.20, 3.23 and 3.33 in \cite{HRS2}. Note that in general one cannot remove the assumption $(A_n) = (A_n)^*$ in assertion $(a)$, whereas $(c)$
holds without this assumption. \hfill \qed \\[3mm]
The notion of a Fredholm sequence was introduced in \cite{RoS7}; see also \cite[Chapter 6]{HRS2}. In the present setting, the Fredholm property of a sequence $(A_n) \in \cS(\eT_u(C))$ means nothing but the invertibility of the coset $(A_n) + \cJ$ in the quotient algebra $\cS(\eT_u(C))/\cJ$, and the results of \cite{RoS7} specify as follows. Let $\sigma_1(a) \le \ldots \le \sigma_n(A) = \|A\|$ denote the singular values of the $n \times n$-matrix $A$.
\begin{theo} \label{t95.42}
Let $(A_n)$ be a sequence in $\cS(\eT_u(C))$ with strong limit $A$. Then \\[1mm]
$(a)$ $(A_n)$ is a Fredholm sequence if and only if $A$ is a Fredholm operator. \\[1mm]
$(b)$ If $A$ is a Fredholm operator and $\dim \ker A = k$, then
\[
\lim_{n \to \infty} \sigma_k(A_n) = 0 \quad \mbox{and} \quad \liminf_{n \to \infty} \sigma_{k+1} (A_n) > 0.
\]
\end{theo}
Assertion $(b)$ allows the numerical determination of the kernel dimension of a Fredholm operator $A \in \eT_u(C)$.
{\small Author's address: \\[3mm]
Steffen Roch, Technische Universit\"at Darmstadt, Fachbereich
Mathematik, Schlossgartenstrasse 7, 64289 Darmstadt,
Germany. \\
E-mail: roch@mathematik.tu-darmstadt.de}
\end{document}